\documentclass{article}

\usepackage{microtype}
\usepackage{graphicx}
\usepackage{subfigure}
\usepackage{booktabs} 
\newcommand{\mat}[1]{{\boldsymbol{{#1}}}} 
\newcommand{\T}[0]{^\top}

\usepackage{hyperref}

\usepackage[accepted]{icml2025}

\makeatletter
\renewcommand{\printAffiliationsAndNotice}[1]{%
\stepcounter{@affiliationcounter}%
{\let\thefootnote\relax\footnotetext{\hspace*{-\footnotesep}\ifdefined\isaccepted #1\fi%
\forloop{@affilnum}{1}{\value{@affilnum} < \value{@affiliationcounter}}{%
\textsuperscript{\arabic{@affilnum}}\csname @affilname\the@affilnum\endcsname}.
\ifdefined\icmlcorrespondingauthor@text
Correspondence to: \icmlcorrespondingauthor@text.
\fi
}}}
\makeatother

\usepackage{amsmath}
\usepackage{amssymb}
\usepackage{mathtools}
\usepackage{amsthm}

\usepackage[capitalize,noabbrev]{cleveref}

\theoremstyle{plain}

\theoremstyle{definition}

\theoremstyle{remark}

\usepackage[textsize=tiny]{todonotes}

\icmltitlerunning{2BSDE Neural Networks for Black--Scholes with Transaction Costs}

\begin{document}

\twocolumn[
\icmltitle{Neural Network Methods in the 2BSDE Formalism for Solving the High-Dimensional Black--Scholes Equation with Transaction Costs}



\icmlsetsymbol{equal}{*}

\begin{icmlauthorlist}
\icmlauthor{Guhesh Kumaran}{yyy}
\end{icmlauthorlist}

\icmlaffiliation{yyy}{Abdus Salam Centre for Theoretical Physics, Imperial College London, Prince Consort Road, London, SW7 2AZ, UK}

\icmlcorrespondingauthor{Guhesh Kumaran}{guhesh.kumaran18@alumni.imperial.ac.uk}


\icmlkeywords{Machine Learning, ICML}

\vskip 0.3in
]



\printAffiliationsAndNotice{}  

\begin{abstract}
A 2BSDE system is written down for a $(D+1)$-dimensional Black--Scholes partial differential equation (PDE) with the \citet{55960847-f1d7-372c-8587-c61145595486} transaction cost model which incorporates the sign of the curvature. The 2BSDE system takes input from the ``speed'' (derivative of the curvature). This system is solved using a single ResNet neural network which learns both the solution to the PDE as well as the degrees of freedom in the curvature matrix. The network is used to learn the price of derivatives such as call options and butterfly spreads on geometric baskets for $D=10,100$ achieving a time-averaged RMSE of less than $0.1$ in each case.
\end{abstract}

\section{Introduction}
\label{sec:Introduction}
Parabolic partial differential equations (PDEs) are ubiquitous in mathematics, physics, and finance. Examples include the heat equation, the Fokker--Planck equation, the Schrödinger equation (after Wick rotation to imaginary time), and the Black--Scholes equation. If the equation is linear then it may be possible to solve it analytically, through Green's function methods etc., but if the equation is nonlinear one may have to turn to numerical methods.

This starts to be challenging for high-dimensional equations which suffer the ``curse of dimensionality'' which refers to the time complexity of lattice methods scaling as the volume of the hypercube. Machine learning methods appear more viable since it is possible to learn PDE solutions and then apply them to new input data without using a lattice.

An alternative representation of a parabolic PDE is the forward-backward stochastic differential equation system (FBSDE) which consists of $D$ forward equations and $1$ backward equation which are generically coupled. The solution to the backward equation is equivalent to the solution of the PDE. This development follows the work of (Feynman-)\citet{92dffd60-7698-3838-9cfd-f3a5e322ac6b}.

Pioneering work by \citet{2017arXiv170604702E} and \citet{2018PNAS..115.8505H} use a neural network to approximate the solution to the backward equation. In \textit{op. cit.} the neural network was used primarily to obtain a point estimate of the PDE solution at some initial point given a terminal condition. \citet{2018arXiv180407010R} then extended this to produce approximations to the PDE solution as a function. 

This paper focuses on the model of \citet{0b9b8115-a8b8-3422-8e1c-a62077de6621} and \citet{b5405564-0209-3283-bc5e-72bd72cfbbed} which gives a parabolic PDE for the price of a derivative on a collection of stocks which obey geometric Brownian motion. The Black--Scholes--Merton model makes some simplifying assumptions in the interest of analyticity, which do not necessarily hold in real markets. This work relaxes the assumption of zero transaction costs when continuously $\Delta$-hedging, which is interesting from both a theoretical and pragmatic point of view. In particular, a simple transaction cost model of \citet{55960847-f1d7-372c-8587-c61145595486} introduces a transaction cost on the change of the hedging size. This introduces a term proportional to the absolute value of the curvature in the Black--Scholes PDE which deforms it to be nonlinear.

An extension of the FBSDE formalism includes a third SDE for the gradient of the PDE solution which contains the curvature. This is called the 2BSDE formalism \cite{2005math......9295C} and it has been used to study similar PDEs in \cite{2019JNS....29.1563B, 2019arXiv190800412P, 2026arXiv260716730Z} for example.

In this work we derive the explicit 2BSDE system for the Black--Scholes equation with Leland's model of transaction costs. This is used to inform the training of a single neural network which learns both the PDE solution and components of its curvature matrix as time varies. This model is then benchmarked on geometric basket call options and butterfly spreads in dimensions $D=10$ and $100$. The nature of the basket essentially reduces the problem to one dimension and hence can be solved with finite-differences. However these products form non-trivial challenges in terms of the change of sign curvature components and in terms of scalability to high dimensions.

\section{2BSDE Systems for Parabolic Partial Differential Equations}
\label{sec:FBSDE}
It was shown, as a consequence of the Feynman--Kac theorem, that solutions to parabolic partial differential equations may be expressed as a coupled system of forward and backward stochastic differential equations \cite{92dffd60-7698-3838-9cfd-f3a5e322ac6b,Pardoux1992BackwardSD, 2005math......9295C,Gobet2020-wh}.

Consider a general form of parabolic PDE with solution $u(t,\vec x): [0,T]\times\mathbb R^D\rightarrow\mathbb R$\begin{align}
    \frac{\partial u}{\partial t} &= a(t,\vec x, u,\nabla u, \nabla^2u) -\vec b(t,\vec x)\cdot\nabla u\nonumber\\&-\frac{1}{2}\operatorname{Tr}\left[\mat c(t,\vec x)\mat c(t,\vec x)\T\nabla^2u\right]\label{eq:genericPDE}
\end{align}where $a$ is a function of $(t,\vec x,u,\nabla u, \nabla^2u)$ but $\vec b, \mat c$ depend only on $(t,\vec x)$. The PDE also has a terminal condition $u(T,\vec x)=f(\vec x)$.

The solution to \cref{eq:genericPDE} can equivalently be represented as a (generically) coupled system of three SDEs called the 2BSDE system \cite{2005math......9295C}
\begin{align}
    d\vec{X}_t &= \vec\mu(t,\vec X_t,Y_t,\vec Z_t)dt +\mat\sigma(t,\vec X_t,Y_t)\cdot d\vec W_t\label{eq:genericfwd}\\
    dY_t &=\varphi(t,\vec X_t,Y_t,\vec Z_t)dt + \vec Z_t\T\cdot \mat\sigma(t,\vec X_t,Y_t)\cdot d\vec W_t\label{eq:genericbwd}\\
    d\vec Z_t&=\vec\psi(t, \vec X_t, Y_t,\vec Z_t, \mat\Gamma_t)dt + \mat\Gamma_t\mat\sigma(t,\vec X_t, Y_t)d\vec W_t\label{eq:genericZ}
\end{align}where $d\vec W_t$ is an increment of (in this case uncorrelated) Brownian motions. The identification of $Y_t=u(t,\vec x),\;\vec Z_t=\nabla_{\vec x}u(t,\vec x),\;\mat\Gamma_t=\nabla^2_{\vec x}u(t,\vec x)$ makes the connection between the PDEs and 2BSDE system explicit.

Using Itô's lemma and the PDE we find the 2BSDE system written in terms of functions that appear in the PDE \cref{eq:genericPDE}\footnote{There is an explicit choice of ``gauge'' which sets $\vec\mu=\vec b$. Other choices include setting $\vec\mu=\vec0$ which turns \cref{eq:genericfwdpde} into a Martingale equation and changes the drift term of \cref{eq:genericbwdpde,eq:genericZpde}.}.
\begin{align}
    d\vec X_t &= \vec bdt + \mat c\cdot d\vec W_t \label{eq:genericfwdpde}\\
    dY_t&=adt+\vec Z_t\T\cdot\mat c\cdot d\vec W_t\label{eq:genericbwdpde}\\
    d\vec Z_t &= \vec\psi dt + \mat\Gamma_t\cdot\mat c\cdot d\vec W_t\label{eq:genericZpde}
\end{align}

and where
\begin{align}
    \vec \psi &= \nabla a+\frac{\partial a}{\partial Y_t}\vec Z_t+\mat\Gamma_t\nabla_{\vec Z_t}a+\operatorname{Tr}_{12}\left[\frac{\partial a}{\partial \mat\Gamma_t}\mat S_t\right]-(\nabla\vec b)\vec Z\nonumber\\&-\frac{1}{2}\operatorname{Tr}_{12}\left[\mat\Gamma_t\nabla(\mat c\mat c\T)\right]
\end{align}with $\operatorname{Tr}_{12}$ indicating a partial trace over the first two indices, $\partial/\partial\mat \Gamma_t=\partial/\partial\Gamma_{ij}$ in components (at time $t$), and $\mat S_t=\nabla^3u$ is the \textit{speed}, borrowing terminology from financial mathematics.

The system has some initial and terminal conditions
\begin{equation}
    \vec X_0=\vec\xi,\; Y_T=f(\vec X_T),\; \vec Z_T=\nabla f(\vec X_T),\;\mat\Gamma_T=\nabla^2f(\vec X_T)
\end{equation}



\section{2BSDE for Black--Scholes Equation and Transaction Costs}

The Black--Scholes (BS) equation is a famous example of a parabolic partial differential equation with particular relevance to quantitative finance. The solution of the BS equation gives the risk-neutral price of a financial derivative, typically options contracts, on some lognormally distributed underlying instruments, typically stocks. In the absence of dividends, the BS equation is written as \begin{equation}
    \frac{\partial u}{\partial t} = ru-r\vec x\cdot\nabla u-\frac{1}{2}\mathrm{Tr}\;\left[\mat\sigma(\vec x)\mat\sigma(\vec x)\T\nabla^2u\right]
\end{equation}
where $t$ is time, $\vec x$ are the prices of the underlyings, $r$ is the risk-free rate, and $\mat\sigma(\vec X)=\mathrm{diag}\;(\sigma_iX_i)$ is the volatility of the underlying. There is a terminal condition which is the payoff of the product.

From here one can identify $a=ru, \vec b = r\vec x$ and so the 2BSDE for the BS equation can be written in components (Latin indices) and suppressing the index for time as\begin{align}
    dX_i &= rX_idt+\sigma_iX_idW_i\label{eq:BSFwd}\\
    dY &=rYdt + \sum_i\sigma_iX_iZ_idW_i\label{eq:BSBwd}\\
    dZ_i &= -\sigma_i^2X_i\Gamma_{ii}dt+\sum_j\Gamma_{ij}\sigma_jX_jdW_j\label{eq:BSZ}
\end{align} where the form of $\vec\psi$ has simplified dramatically and the speed is not needed.

\subsection{Transaction Costs}
\label{sec:transaction}
The BS model ignores transaction costs which can become significant when continuously hedging. The transaction cost model of \citet{55960847-f1d7-372c-8587-c61145595486} offers a simple deformation of the BS model which is easy to implement in the 2BSDE system. For alternative transaction costs models \citep[see e.g.][]{Barles1998-sz, Lepinette2012-ae, Herdegen2021-nf, Saeed2026-cg}.

The $D$ assets follow independent geometric Brownian motions as in \cref{eq:BSFwd} with re-hedging occurring at time intervals $\Delta t_H$. The transaction cost is paid on the amount the hedged position moves by, which is $\sum_j\Gamma_{ij}\Delta X_j+O(\Delta t_H)=\sum_{j}\Gamma_{ij}\sigma_jX_j\Delta W_j+O(\Delta t_H)$ units where $\Delta W_j\sim N(0,\Delta t_H)$.

The total transaction cost is given by
\begin{equation}
  \mathrm{TC}=\sum_{i=1}^{D}\frac{k_i}{2}\,X_i\,\lvert\Delta Z_i\rvert,
\end{equation}
where $k_i$ is paid for a ``round-trip'' per unit value traded. The expected transaction cost per unit (hedging) time is
\begin{equation}
  \frac{\mathbb{E}[\mathrm{TC}]}{\Delta t_H}
  =\frac12\sum_{i=1}^{D}\mathrm{Le}_i\,\sigma_i X_i
   \Big(\sum_{j=1}^{D}\sigma_j^2X_j^2\,\Gamma_{ij}^2\Big)^{1/2}
\end{equation} with the definition of the Leland number \cite{Marco01121994}
\begin{equation}
    \mathrm{Le}_i=\sqrt{\frac{2}{\pi\Delta t_H}}\frac{k_i}{\sigma_i} .
\end{equation}

A further simplification is made by keeping only the diagonal contributions which simplifies the expected transaction cost per unit time to
\begin{equation}
\frac{\mathbb E[\mathrm{TC}]}{\Delta t_H}=\frac12\sum_i\mathrm{Le}_i\,\sigma_i^2X_i^2\,|\Gamma_{ii}|
\end{equation}
This results in a drastic simplification of the 2BSDE system. In the BS PDE the transaction cost term amounts to replacing each $\sigma_i^2$ in the BS equation by $\hat\sigma_i^2=\sigma_i^2\big(1+\mathrm{Le}_i\,\mathrm{sign}\,\Gamma_{ii}\big)$ where the sign of $\Gamma_{ii}$ makes the equation nonlinear. However it is more convenient to separate out these terms as
\begin{align}
\frac{\partial u}{\partial t}&=ru-\frac12\sum_i\mathrm{Le}_i\sigma_i^2X_i^2|\Gamma_{ii}|-r\sum_jX_j\frac{\partial u}{\partial X_j}\nonumber\\
&-\frac12\sum_j\sigma_j^2X_j^2\Gamma_{jj}.
\end{align}which then takes the form of a Black--Scholes--Barenblatt equation similar to the one studied in \cite{2019JNS....29.1563B} where $\sigma^2_{i,\max,\min}=\sigma_i^2(1\pm\operatorname{Le}_i)$.

In the notation of the general system, $a=rY-\frac12\sum_i\mathrm{Le}_i\sigma_i^2X_i^2|\Gamma_{ii}|$, $b_i=rX_i$ and $\mat c=\mathrm{diag}(\sigma_iX_i)$, so the 2BSDE system is
\begin{align}
dX_i&=rX_i\,dt+\sigma_iX_i\,dW_i,\label{eq:LelandFwd}\\
dY&=\Big(rY-\frac12\sum_i\mathrm{Le}_i\sigma_i^2X_i^2|\Gamma_{ii}|\Big)dt\nonumber\\
&\quad+\sum_iZ_i\sigma_iX_i\,dW_i,\label{eq:LelandBwd}\\
dZ_i&=-\Big[\sigma_i^2X_i\big(1+\mathrm{Le}_i\,\mathrm{sign}\,\Gamma_{ii}\big)\Gamma_{ii}\nonumber\\
&\quad+\frac12\sum_j\mathrm{Le}_j\sigma_j^2X_j^2\,\mathrm{sign}(\Gamma_{jj})\,S_{jji}\Big]dt\nonumber\\
&\quad+\sum_j\Gamma_{ij}\sigma_jX_j\,dW_j,\label{eq:LelandZ}
\end{align}
where the speed now explicitly appears in the 2BSDE system.
\section{Deep Neural Networks to solve PDEs}
\label{sec:DNNtosolvePDE}
Having developed the 2BSDE system in \cref{eq:LelandFwd,eq:LelandBwd,eq:LelandZ} it remains to solve them. It is clear that the forward equation \eqref{eq:LelandFwd} is uncoupled and so it may be solved by Euler-Maruyama discretisation. The backward equation \eqref{eq:LelandBwd} and the equation for $Z_i$ \eqref{eq:LelandZ} must be adapted to the filtration induced by the forward equation which complicates solving them.

One approach is to use a neural network to approximate the solution $Y(t,\vec X_t)=u(t,\vec X_t)\simeq u^\theta(t,\vec X_t)$ where $\theta$ are a collection of parameters in the neural network which are shared across all time steps. Then $\vec Z_t$ may be computed with automatic differentiation. In principle one can compute $\mat \Gamma_t$ and $\mat S_t$ through further automatic differentiation but this tends to be computationally expensive in large dimensions. Instead \cref{eq:LelandZ} can be used to learn $\mat\Gamma^\theta(t,\vec X_t)\simeq\mat\Gamma(t,\vec X_t)$ and compute $\mat S_t(t,\vec X_t)=\nabla_{\vec X}\mat\Gamma^\theta(t,\vec X_t)$ through its automatic differentiation. Note that the same neural network parameters $\theta$ are used to learn both $u^\theta(t,\vec X_t)$ and $\mat\Gamma^\theta(t,\vec X_t)$ as opposed to two separate neural networks that learn each individually.

Further let us discretise the time interval $[0,T]$ into $N$ parts, $0=t_0<t_1<\dots<t_N=T$, with steps of $\Delta t_n=t_{n+1}-t_n$ and
$\Delta \vec W_n=\vec W_{t_{n+1}}-\vec W_{t_n}$. Along each simulated path,
\begin{equation}
  (Y_n,\,Z_n,\,\Gamma_n)=\big(u^\theta,\,\nabla_X u^\theta,\,\Gamma^\theta\big)(t_n,X_n).
\end{equation}
A subscript $n$ on a coefficient denotes evaluation at time $t_n$.

The neural network needs a loss function which is partly taken from minimising the $L^2$ variation in time steps in $u^\theta(t,\vec X_t)$ and in $\nabla_{\vec X}u^\theta(t,\vec X_t)$ \cite{2018arXiv180407010R}. It is sufficient to use Euler-Maruyama discretisation
\begin{align}
  \vec X_{n+1} &= \vec X_n+\vec \mu_n\,\Delta t_n+\mat\sigma_n\,\Delta \vec W_n,\\
  \widetilde Y_{n+1} &= Y_n+\varphi_n\,\Delta t_n+\vec Z_n\T\mat\sigma_n\,\Delta \vec W_n,\\
  \widetilde{ \vec Z_{n+1}} &= \vec Z_n+\vec\psi_n\,\Delta t_n+\mat\Gamma_n\,\mat\sigma_n\,\Delta \vec W_n,
\end{align}
The loss also has a ``shooting term'' which seeks to minimise $L^2$ deviation in $u^\theta(0,\vec X_0)$ to the ``true'' value \cite{2017arXiv170604702E}. This requires an estimate of the terminal value $\mathcal Y_N$ from a given $u^\theta(0,\vec X_0)=Y_0$ by adding the incremental values at each time step as
\begin{equation}
  \mathcal Y_N = Y_0+\sum_{n=0}^{N-1}\big(\varphi_n\,\Delta t_n+Z_n^{\top}\sigma_n\,\Delta W_n\big).
\end{equation}and then penalising the $L^2$ deviation of $\mathcal Y_N$ from $f(\vec X_N)$.

Finally, there are terms in the loss which seek to minimise the $L^2$ deviation in the terminal conditions for $Y_N, \vec Z_N,\mat\Gamma_N$. In summary, the loss over a batch of $M$ paths, indexed by $m$, is
\begin{equation}
  \mathcal{L}(\theta)=\sum_{m=1}^{M}\Big(\sum_{n=0}^{N-1}\ell^m_n+\ell^m_N\Big),\label{eq:loss}
\end{equation}
with step and terminal terms
\begin{align}
  \ell^m_n &= \big|Y^m_{n+1}-\tilde Y^m_{n+1}\big|^2
             +\lambda_Z\,\big\|Z^m_{n+1}-\tilde Z^m_{n+1}\big\|^2,\\
  \ell^m_N &= \big|Y^m_N-f(X^m_N)\big|^2+\big|\mathcal Y^m_N-f(X^m_N)\big|^2
             \nonumber\\
           &\hphantom{{}={}}+\big\|Z^m_N-\nabla f(X^m_N)\big\|^2
             \nonumber\\
           &\hphantom{{}={}}+\lambda_\Gamma\,\big\|\Gamma^m_N-\nabla^2 f(X^m_N)\big\|_F^2,
\end{align}
where $\lambda_Z=\lambda_\Gamma=1$ if the $\Gamma$-head is switched on, which is the case for most experiments in this work, and zero otherwise.

\section{Neural Network Architecture and Training}
The choice of neural network architecture made here is the ResNet architecture \cite{he2016deep} with a sine activation function \cite{sitzmann2020implicitneuralrepresentationsperiodic}. The sine activation function is particularly useful to ensure smooth derivatives at any order.

Crucially, in the loss function we do not allow gradients of $\mat\Gamma^\theta_n$ to flow through $\varphi_n$ since this tends to bias $\mat\Gamma^\theta_n$.

The Adam optimiser \cite{Kingma2014AdamAM} was used with a two stage learning schedule. The first stage occurs with a learning rate which is determined by a learning rate range test \cite{smith2017cyclical} which can be implemented in the \texttt{fastai} package \cite{info11020108}. The chosen learning rate $\eta^\ast$ is the ``valley'' learning rate. The first stage occurs for a maximum of 40000 iterations with an early stopping and the second stage occurs with a learning rate of $\eta^\ast/10$ with a maximum of 20000 iterations with the same early stopping condition. The early stopping condition used here comes from tracking the mean training loss every 100 iterations and stopping once the median of the last ten of these mean losses has not improved by 2\% in the last 2000 iterations. 

Since both $u^\theta(t,\vec X_t)$ and $\mat\Gamma^\theta(t,\vec X_t)$ are learnt from the same neural network there are some relative scales between their values which are important to mind in the loss function. This is to be mitigated by tuning the coefficients $\lambda_Z$ and $\lambda_\Gamma$ which are both set to one here. 


The training occurs with $N=50$ time steps in batch sizes of $M=100$.

The code used in this work is available at \url{https://github.com/guheshk/Leland2BSDE} and extends the code base of \citet{2018arXiv180407010R} and \citet{2019arXiv191011623G}.

Further details about the network, its training, and hyperparameters can be found in \cref{app:hyper}.

\section{Results}
Many cases in different dimensions are studied here to show both agreement with the reference solutions and also the ability of a neural network to learn the effects of transaction costs.

The experiments are run using some common values of a constant risk-free rate of $r=0.1$, a maturity of $T=1$, a constant re-hedging interval of $\Delta t_H=1/52$. Some parameters are taken to be identical in the $D$ dimensional cases for simplicity. In particular we take $\forall i\in\{1,\cdots,D\},\;X_0^i=100,\;\sigma_i=0.2,\;k_i=0.005,$ which implies $\textrm{Le}_i=0.1438$. We also assume no correlations between Brownian motions.

The geometric basket $G_t$ is defined as \begin{equation}
    G_t=\prod_{i=1}^DX_i^{w_i},\quad \sum_{i=1}^Dw_i=1, \quad w_i\geq 0
\end{equation}where in this work we take $w_i=1/D$. The geometric basket reduces problems in high dimensions to a problem in one dimension which can be solved by finite differences which makes it a suitable reference for the neural network learnt solutions. 

The intrinsic nature of the geometric basket which grows as the geometric mean of the constituent assets leads to changing hedging ratios for each asset within the basket to maintain geometric growth. These incur transaction costs and make the geometric basket suitable for the study in this work. There is a quantity $\mathcal Q$ defined as \begin{equation}
    \mathcal Q=G^2\frac{\partial^2 V}{\partial G^2}-(D-1)G\frac{\partial V}{\partial G},\quad V(t,G_t)=u(t,\vec X_t)
\end{equation}that essentially replaces the role of $\mat\Gamma_t$ in the Black--Scholes equation in $D$ underlyings with the Black--Scholes equation in the basket level only. The Leland term is dependent on the sign of $\mathcal Q$ which for $D>1$ can be positive or negative, making geometric baskets call options non-trivial tests of the Leland model in higher dimensions. Further information about geometric baskets is provided in \cref{sec:GeomBasket}. This one dimensional Black--Scholes equation can be solved by the method of finite-differences and will serve as the reference solution for the neural network. Further details about the finite-difference solver are in \cref{app:fd}.

The results are summarised in \cref{tab:summary} and we will refer to it throughout this section.
\begin{table*}[t]
\caption{Results of the experiments conducted in this work on equally weighted geometric basket calls and butterfly spreads. The reference values denoted with superscript ``ref'' are computed from the finite-difference solution of the one-dimensional Leland PDE. The ``final'' label refers to the value evaluated by the network on its final iteration of training and ``mean'' refers to the mean value $\pm$ sample standard deviation of the last 20 checkpoints (2000 iterations). The RMSE of $Y_t$ is evaluated over 100 test paths at 51 time points. The proportion of the friction premium captured by the model is computed as $(Y_0-Y_0^{\mathrm{fr}})/(Y_0^{\mathrm{ref}}-Y_0^{\mathrm{fr}})$ with frictionless prices $Y_0^{\mathrm{fr}}$ coming from finite difference methods without friction. The total training time (on a standard CPU) is presented along with the time for 100 iterations in brackets. There is a special run denoted by $^\dagger$ which was trained without stopping gradients of $\Gamma$ flowing through the loss term for deviations in $Y_n^m$.}
\label{tab:summary}
\vskip 0.1in
\begin{center}
\begin{small}
\begin{tabular*}{\textwidth}{@{\extracolsep{\fill}}lcccccccc@{}}
\toprule
 & & & \multicolumn{2}{c}{$Y_0$ (NN)} & Rel.\ error & & Premium & Time \\
\cmidrule(lr){4-5}
Product & $D$ & $Y_0^{\mathrm{ref}}$ & final & mean & final\,/\,mean & RMSE & final\,/\,mean & (s/100 it.) \\
\midrule
Call & 10 & 8.1954 & 8.1533 & $8.179\pm0.023$ & $-0.51\%$\,/\,$-0.20\%$ & 0.057 & 78\%\,/\,91\% & 2.5 h (35) \\
Butterfly & 10 & 2.9065 & 2.8835 & $2.894\pm0.028$ & $-0.79\%$\,/\,$-0.44\%$ & 0.090 & 91\%\,/\,95\% & 2.0 h (36) \\
Butterfly (no $\Gamma$ head) & 10 & 2.9065 & 2.8703 & $2.895\pm0.029$ & $-1.24\%$\,/\,$-0.40\%$ & 0.058 & 85\%\,/\,95\% & 15.0 h (195) \\
\midrule
Call & 100 & 7.8354 & 7.8245 & $7.774\pm0.033$ & $-0.14\%$\,/\,$-0.78\%$ & 0.068 & 96\%\,/\,78\% & 11.4 h (149) \\
Butterfly & 100 & 2.0078 & 1.9757 & $1.969\pm0.025$ & $-1.60\%$\,/\,$-1.93\%$ & 0.049 & 88\%\,/\,85\% & 12.9 h (144) \\
\midrule
Butterfly$^\dagger$ & 10 & 2.9065 & 3.1828 & $3.190\pm0.030$ & $+9.51\%$\,/\,$+9.75\%$ & 0.297 & 213\%\,/\,216\% & 2.9 h (40) \\
\bottomrule
\end{tabular*}
\end{small}
\end{center}
\vskip -0.1in
\end{table*}
\begin{figure*}[h]
\centering
\subfigure[$D=10$ call]{\includegraphics[width=0.32\textwidth]{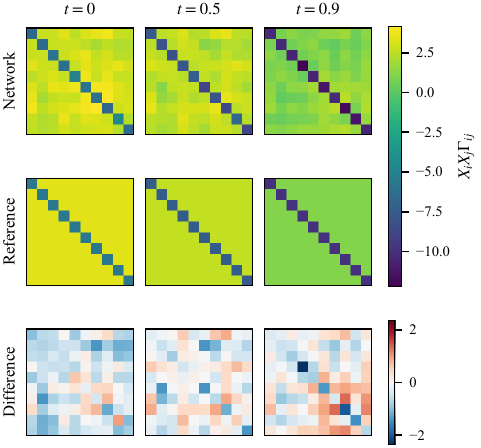}\label{fig:MapCallD10}}\hfill
\subfigure[$D=10$ butterfly]{\includegraphics[width=0.32\textwidth]{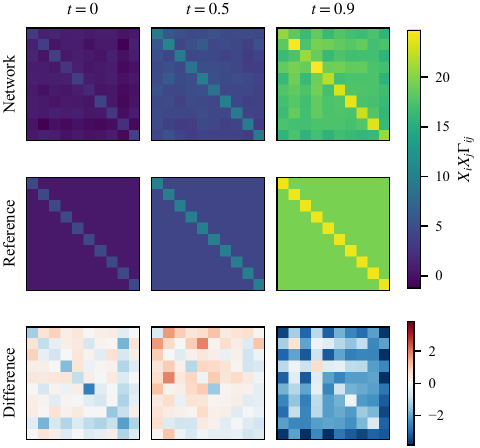}\label{fig:MapFlyD10}}\hfill
\subfigure[$D=100$ butterfly]{\includegraphics[width=0.32\textwidth]{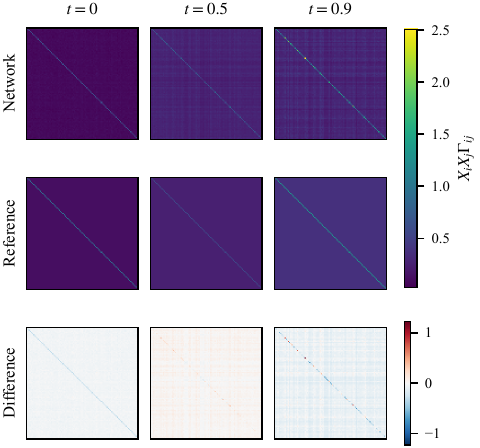}\label{fig:MapFlyD100}}
\caption{Cash gamma $X_iX_j\Gamma_{ij}$ on one test path at $t=0, 0.5, 0.9$ for the $D=10$ geometric basket call and the $D=10, 100$ geometric basket butterfly spreads. The $\Gamma$-head output (top), finite difference reference (middle) and their difference (bottom).}
\label{fig:HeatMaps}
\end{figure*}

\subsection{The $D=10$ geometric basket}
\label{sec:D10Basket}
\begin{figure*}[h]
\centering
\subfigure[Sample paths]{\includegraphics[width=0.32\textwidth]{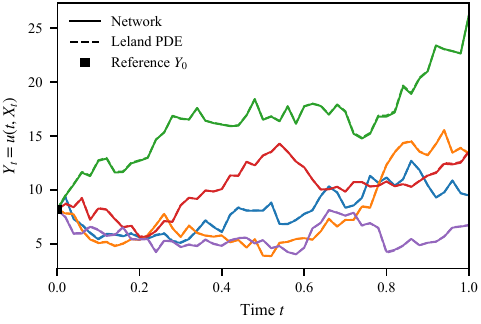}\label{fig:BasketD10Paths}}\hfill
\subfigure[RMSE over time]{\includegraphics[width=0.32\textwidth]{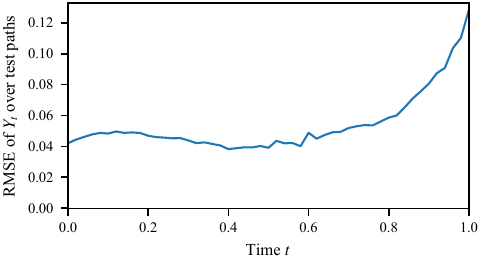}\label{fig:BasketD10RMSE}}\hfill
\subfigure[Error in $X_iX_j\Gamma_{ij}$]{\includegraphics[width=0.32\textwidth]{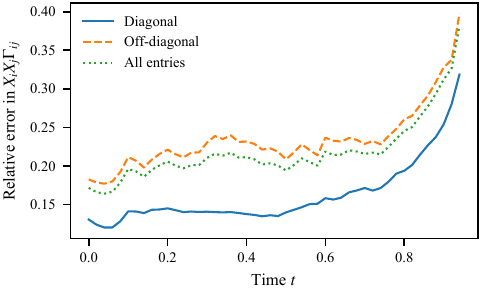}\label{fig:BasketD10GammaError}}
\caption{$D=10$ geometric basket call. (a) Leland PDE solution (dashed) and network (solid) on sample paths. (b) RMSE of $Y_t$ over the test paths at each time. (c) Relative error of the learned cash gamma against the exact one at each time, for the diagonal, off-diagonal and all entries.}
\label{fig:BasketD10}
\end{figure*}


The first product studied here is the $D=10$ geometric basket which has a payoff given by \begin{equation}
    f(G_T)=(G_T-K)^+
\end{equation}where $K=100$ is the strike price. 

The solution for the option price can be solved by a finite difference solver giving $Y_0^{\mathrm{ref}}=8.1954$, compared with the frictionless price of $8.0057$. The $\Gamma$-head outputs the $D(D+1)/2=55$ degrees of freedom of the symmetric matrix $\mat\Gamma^\theta$, the ResNet architecture has 214840 learnable parameters and the learning rate range test finds $\eta^*=1.51\times10^{-4}$.

The network trained for a total of 25900 iterations which took 2.5 hours on a standard CPU. The final option price $Y_0$ has a $-0.51\%$ relative error, and averaged over the last 2000 iterations gives $Y_0=8.179$ with a $-0.20\%$ relative error. The pathwise RMSE along sample trajectories averages to 0.057. Figure \ref{fig:BasketD10Paths} shows that the network and the Leland PDE solution are almost indistinguishable along the sample paths. The average RMSE over the test paths is shown in Figure \ref{fig:BasketD10RMSE}, where it stays at around $0.04$--$0.05$ until $t\approx0.6$ and then increases steadily to around $0.13$ at maturity.


The numerical values of $\Gamma_{ij}$ tend to be rather small in high dimensions and so a scaling known as the ``cash gamma'' $C_{ij}=X_iX_j\Gamma_{ij}$ is used here. Figure \ref{fig:MapCallD10} shows the learnt and exact cash gamma on one test path at $t=0$, $0.5$ and $0.9$. The network learnt (cash) gamma agrees closely with the exact (cash) gamma at each of the three times. In particular the negative diagonal terms and positive off-diagonal terms are seen in the network. The heatmap of the differences remain quite small at early times with some larger differences at $t=0.9$ which is close to expiry. The relative errors in Figure \ref{fig:BasketD10GammaError} show that the diagonal entries, which are the ones entering the transaction-cost term, are learnt to within about $14\%$ up to $t\approx0.5$, rising to around $30\%$ near maturity. The off-diagonal entries have a larger error of around $20\%$, rising to around $40\%$ near maturity.
\subsection{$D=10$ Butterfly Spread}
\label{sec:D10Butterfly}

\begin{figure*}[h!]
\centering
\subfigure[Sample paths]{\includegraphics[width=0.32\textwidth]{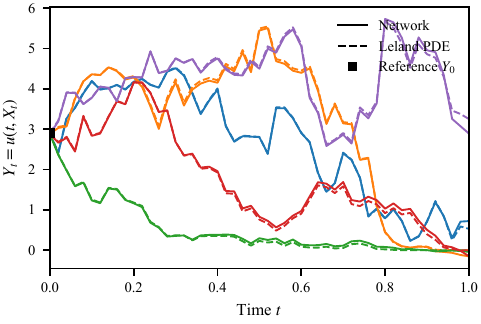}\label{fig:ButterflyD10Paths}}\hfill
\subfigure[RMSE over time]{\includegraphics[width=0.32\textwidth]{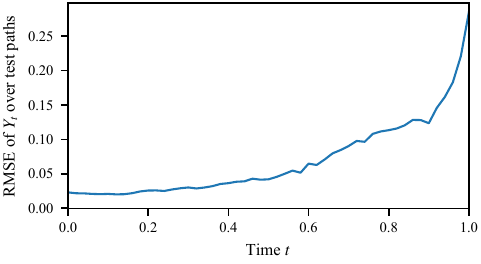}\label{fig:ButterflyD10RMSE}}\hfill
\subfigure[Error in $X_iX_j\Gamma_{ij}$]{\includegraphics[width=0.32\textwidth]{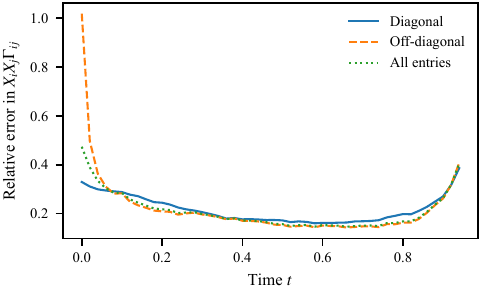}\label{fig:ButterflyD10GammaError}}
\caption{$D=10$ geometric basket butterfly spread (strikes 90/100/110). (a) Leland PDE solution (dashed) and network (solid) on sample paths. (b) RMSE of $Y_t$ over the test paths at each time. (c) Relative error of the learned cash gamma against the exact one at each time, for the diagonal, off-diagonal and all entries.}
\label{fig:ButterflyD10}
\end{figure*}

A more demanding example is the $D=10$ geometric basket butterfly spread. The payoff takes the following form \begin{equation}
    f(G_T)=(G_T-K_1)^+-2(G_T-K_2)^++(G_T-K_3)^+
\end{equation}where $K_{1,2,3}=90,100,110$. The basket is defined on the same assets and weightings as in \cref{sec:D10Basket}. This time there are two sources which affect the change of sign of $\mat\Gamma_t$; there is the intrinsic change of sign in the basket and also the change in sign associated with the payoff structure. The geometric basket reference price follows the one-dimensional finite difference solver, which gives $Y_0^{\mathrm{ref}}=2.9065$. This is higher than the frictionless price of $2.6629$.

The network trained for a total of 19700 iterations which took 2.0 hours on a standard CPU. The final option price $Y_0$ has a $0.79\%$ relative error whereas averaging it over the last 2000 iterations gives $Y_0=2.894$ with a $0.44\%$ relative error. Agreement of the NN paths with the PDE paths for the whole solution is strong, indicated by the average RMSE of 0.090 on the sample paths. Figure \ref{fig:ButterflyD10Paths} visually shows this agreement except close to maturity where there are some deviations. The average RMSE over the test paths is shown in Figure \ref{fig:ButterflyD10RMSE}, where it grows steadily from around $0.02$ at $t=0$ up to $0.13$ around $t=0.9$ and then sharply increases for $t>0.9$. This is most likely due to the terminal conditions fixing the end points of trajectories and resulting in sharp jumps in the learnt solution to satisfy them.


The learnt and exact cash gamma on one test path are shown in Figure \ref{fig:MapFlyD10}. The network reproduces the structure of the exact matrix at all three times, including the growth of the cash gamma towards maturity on this path, although at $t=0.9$ the learnt entries are slightly smaller than the exact ones. The relative errors in Figure \ref{fig:ButterflyD10GammaError} show that the diagonal entries are learnt to within about $30\%$ at $t=0$, falling to around $17\%$ between $t\approx0.4$ and $t\approx0.8$ and rising to around $40\%$ near maturity. The off-diagonal error is large at $t=0$, where the exact off-diagonal entries are close to zero, and is then similar to the diagonal error.

\subsection{$D=10$ butterfly spread with and without the $\Gamma$-head}
\label{sec:ButterflyD10NoHead}
\begin{figure}[h]
\centering
\includegraphics[width=\linewidth]{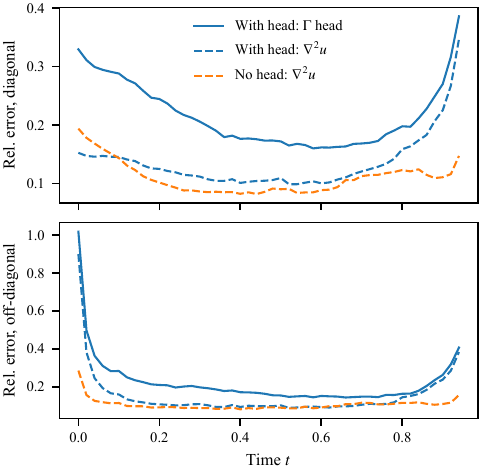}
\caption{Relative error of the cash gamma over the test paths at each time for the $D=10$ basket butterfly spread, for the diagonal (top) and off-diagonal (bottom) entries, trained with a $\Gamma$ head (solid: head output; dashed: the network's Hessian $\nabla^2u$) and without one (dashed: $\nabla^2u$).}
\label{fig:ButterflyD10HeadVsNoHead}
\end{figure}

In this section we study the same $D=10$ butterfly spread but turn off the $\Gamma$-head and set $\lambda_Z=\lambda_\Gamma=0$ in the training loss. Instead of learning $\mat\Gamma_t$ we use the automatic differentiation of the learned solution. The learning rate range test finds $\eta^*=2.16\times10^{-4}$.

A striking difference is that in this case training took 15.0 hours (at 195s per 100 iterations) on a standard CPU compared to 2.0 hours (at 36s per 100 iterations) when the $\Gamma$-head was used. The relative error in $Y_0$ is worse in this case at $-1.24\%$ on the final iteration and comparable at $-0.40\%$ in the average over the final 2000 iterations compared with the $\Gamma$-head. The pathwise RMSE is lower without the $\Gamma$-head at $0.058$ and higher with it at $0.090$.

Plots of the relative error in the diagonal and off-diagonal elements of the cash $\mat\Gamma_t$ are shown in \cref{fig:ButterflyD10HeadVsNoHead}. It is apparent that the $\Gamma$-head has higher relative errors in both the diagonal and off-diagonal elements than compared to using the automatic differentiation of the network solution (for both the cases with and without the $\Gamma$-head).

\subsection{$D=100$ geometric basket call}
\label{sec:BasketD100}

\begin{figure*}[h!]
\centering
\subfigure[Sample paths]{\includegraphics[width=0.32\textwidth]{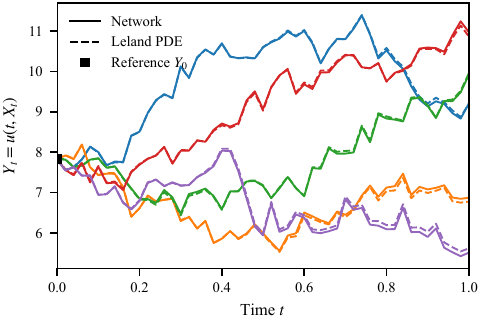}\label{fig:BasketD100Paths}}\hfill
\subfigure[RMSE over time]{\includegraphics[width=0.32\textwidth]{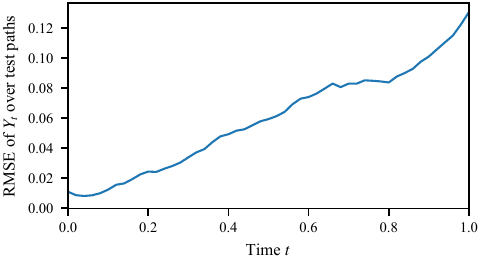}\label{fig:BasketD100RMSE}}\hfill
\subfigure[Error in $X_iX_j\Gamma_{ij}$]{\includegraphics[width=0.32\textwidth]{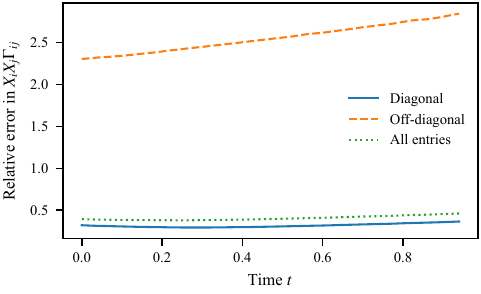}\label{fig:BasketD100GammaError}}
\caption{$D=100$ geometric basket call. (a) Leland PDE solution (dashed) and network (solid) on sample paths. (b) RMSE of $Y_t$ over the test paths at each time. (c) Relative error of the learned cash gamma against the exact one at each time, for the diagonal, off-diagonal and all entries.}
\label{fig:BasketD100}
\end{figure*}

The method is next scaled up to the $D=100$ geometric basket call option, with the same strike $K=G_0=100$ as before. The basket volatility is now $\sigma_G=\sigma/\sqrt{D}=0.02$ and the reference solution from the finite difference solver gives $Y_0^{\mathrm{ref}}=7.8354$, compared with the frictionless price of $7.5558$. The $\Gamma$-head now outputs the $D(D+1)/2=5050$ independent entries of $\Gamma$ and the ResNet architecture has 1521595 learnable parameters.


The network trained for a total of 27600 iterations which took 11.4 hours on a standard CPU. The final option price achieved by the network is $Y_0=7.825$ which has a $0.14\%$ relative error, however when averaged over the last 2000 iterations, $Y_0=7.774$ with a $0.78\%$ relative error. The pathwise RMSE along sample trajectories averages to 0.068. Figure \ref{fig:BasketD100Paths} shows that the network and the finite difference solution are almost indistinguishable along the sample paths with only minor deviations near maturity. This is corroborated by the plot of the average RMSE against time shown in \cref{fig:BasketD100RMSE} which starts at around $0.01$ at $t=0$ and then increases almost linearly to around $0.13$ at maturity.


The plots of the relative errors in the cash gamma components at different times are given in Figure \ref{fig:BasketD100GammaError}. The relative errors in the diagonal terms are learnt within $40\%$ at all times. However the relative errors of the off-diagonal terms are markedly higher starting at around $2.3$ at $t=0$ and increasing steadily up to around $2.9$ at maturity. This is mainly because the off-diagonal terms are approximately zero and so any fluctuation about zero results in significant relative errors. There is no clipping of small values implemented here since the off-diagonal terms are not so important for the transaction costs.

\subsection{$D=100$ geometric basket butterfly spread}
\label{sec:ButterflyD100}

\begin{figure*}[t!]
\centering
\subfigure[Sample paths]{\includegraphics[width=0.32\textwidth]{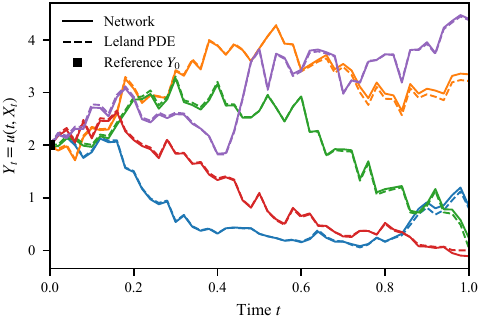}\label{fig:ButterflyD100Paths}}\hfill
\subfigure[RMSE over time]{\includegraphics[width=0.32\textwidth]{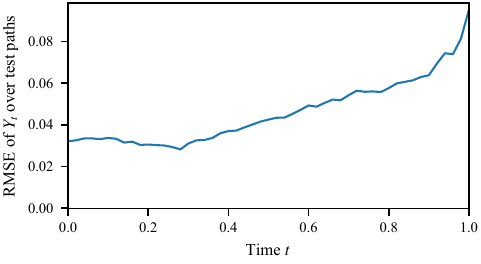}\label{fig:ButterflyD100RMSE}}\hfill
\subfigure[Error in $X_iX_j\Gamma_{ij}$]{\includegraphics[width=0.32\textwidth]{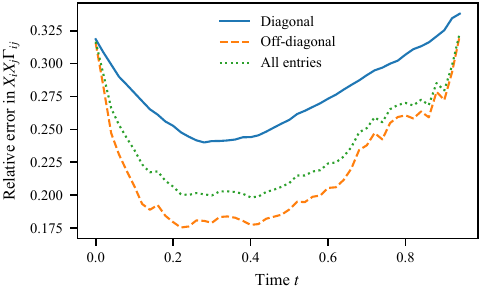}\label{fig:ButterflyD100GammaError}}
\caption{$D=100$ geometric basket butterfly spread (strikes 90/100/110). (a) Leland PDE solution (dashed) and network (solid) on sample paths. (b) RMSE of $Y_t$ over the test paths at each time. (c) Relative error of the learned cash gamma against the exact one at each time, for the diagonal, off-diagonal and all entries.}
\label{fig:ButterflyD100}
\end{figure*}

The final example is the $D=100$ geometric basket butterfly spread on the same basket with the same payoff as used in $D=10$. The one-dimensional finite difference solver gives $Y_0^{\mathrm{ref}}=2.0078$, compared with the frictionless price of $1.7493$.

The network trained for a total of 32200 iterations which took 12.9 hours on a standard CPU. The final option price achieved by the network is $Y_0=1.976$ corresponding to a $1.60\%$ relative error. When averaged over the last 2000 iterations the network achieves $Y_0=1.969$ which is a $1.93\%$ relative error. The pathwise RMSE along sample trajectories averages to 0.049 which is the lowest of all examples studied here. This is seen clearly in plots of some sample paths in Figure \ref{fig:ButterflyD100Paths} where once again the learnt solution and finite difference solution are visually indistinguishable. A plot of the average RMSE over the test paths is given in \cref{fig:ButterflyD100RMSE} where the RMSE starts at around $0.03$ for $0<t<0.3$ and then increases linearly to around $0.06$ until $t=0.9$ with a final sharp increase to around $0.10$ at maturity.


The network is able to learn the structure of the cash gamma at all three times as shown in \cref{fig:MapFlyD100}. The diagonal consists of some larger positive values with the off-diagonal remaining positive but slightly smaller. In general there are only slight differences in the values learnt and the exact values, but this difference does increase near maturity. On average, the relative errors along time, shown in Figure \ref{fig:ButterflyD100GammaError}, also show this trend. The relative errors at $t=0$ and $t=1$ are at roughly $0.325$ for both diagonal and off-diagonal elements. The relative error decreases until around $t=0.4$, with a larger decrease for off-diagonal elements than the diagonal ones. The lowest relative error for the diagonal elements is around $0.25$ and around $0.175$ for the off-diagonal ones.

In this case it would not be feasible to run an experiment comparing the $\Gamma$-head network to one which uses automatic differentiation only given that in the $D=10$ case not using the $\Gamma$-head was already five times slower.

\section{Discussion}
In this work a 2BSDE system is studied for a Black--Scholes model with a deformation that incorporates a transaction cost on changing hedging ratios following the model of \citet{55960847-f1d7-372c-8587-c61145595486}. The transaction cost appears in the Black--Scholes PDE through the sign of the curvature and additionally the 2BSDE system requires the ``speed''.

A single neural network method is proposed which learns both the PDE solution and the degrees of freedom in its curvature matrix. This follows ideas and code from \cite{2005math......9295C, 2017arXiv170604702E, 2018PNAS..115.8505H, 2018arXiv180407010R, 2019JNS....29.1563B,2019arXiv191011623G}.

The neural network is tested on a range of  examples including a geometric basket call option and a butterfly spread on a geometric basket in dimensions $D=10,100$. These intrinsically include a change in sign $\mat\Gamma_t$ which is relevant for this model. In each case there is strong agreement with the learnt PDE solution, the relative error in $Y_0$ was between $0.2-1.9\%$ and the RMSE was between $0.05-0.09$. The components of $\mat\Gamma_t$ are learnt sufficiently well to solve the PDE which is the primary goal. The inclusion of the $\Gamma$-head drastically reduces computation time for high-dimensional examples, positioning this method for scalability in high dimensions.

This work examines only simple neural network architectures and somewhat simple derivative payoffs. The examples in this work also use the same coefficients for the forward system and baskets. So future work could entail optimising the methods employed here for specific products traded in real life such as exotic options or arithmetic basket options, which do not have closed form solutions or reduce to a lower dimensional problem. One avenue for studying other types of product could include generating solutions by repeated Monte Carlo simulations at different time steps. Further improvements in the solvers could be found by using different neural network architectures or through hyperparameter optimisation of learning rates or the coefficients $\lambda_Z,\lambda_\Gamma$ that appear in the loss function.
\section*{Acknowledgements}
The author thanks Panos Parpas for encouragement in preparing this preprint.





\bibliography{references}

@BOOK{Gobet2020-wh,
  title     = "{Monte-Carlo} methods and stochastic processes",
  author    = "Gobet, Emmanuel",
  publisher = "CRC Press",
  month     =  sep,
  year      =  2020,
  address   = "London, England",
  language  = "en"
}

@ARTICLE{2018arXiv180407010R,
       author = {{Raissi}, Maziar},
        title = "{Forward-Backward Stochastic Neural Networks: Deep Learning of High-dimensional Partial Differential Equations}",
      journal = {arXiv e-prints},
         year = 2018,
        month = apr,
          eid = {arXiv:1804.07010},
        pages = {arXiv:1804.07010},
          doi = {10.48550/arXiv.1804.07010},
archivePrefix = {arXiv},
       eprint = {1804.07010},
 primaryClass = {stat.ML},
       adsurl = {https://ui.adsabs.harvard.edu/abs/2018arXiv180407010R}
}

@ARTICLE{2017arXiv170604702E,
       author = {{E}, Weinan and {Han}, Jiequn and {Jentzen}, Arnulf},
        title = "{Deep learning-based numerical methods for high-dimensional parabolic partial differential equations and backward stochastic differential equations}",
      journal = {arXiv e-prints},
         year = 2017,
        month = jun,
          eid = {arXiv:1706.04702},
        pages = {arXiv:1706.04702},
          doi = {10.48550/arXiv.1706.04702},
archivePrefix = {arXiv},
       eprint = {1706.04702},
 primaryClass = {cs.NA},
       adsurl = {https://ui.adsabs.harvard.edu/abs/2017arXiv170604702E}
}

@ARTICLE{2018PNAS..115.8505H,
       author = {{Han}, Jiequn and {Jentzen}, Arnulf and {E}, Weinan},
        title = "{Solving high-dimensional partial differential equations using deep learning}",
      journal = {Proceedings of the National Academy of Science},
         year = 2018,
        month = aug,
       volume = {115},
       number = {34},
        pages = {8505-8510},
          doi = {10.1073/pnas.1718942115},
archivePrefix = {arXiv},
       eprint = {1707.02568},
 primaryClass = {cs.NA},
       adsurl = {https://ui.adsabs.harvard.edu/abs/2018PNAS..115.8505H}
}

@ARTICLE{2005math......9295C,
       author = {{Cheridito}, Patrick and {Mete Soner}, H. and {Touzi}, Nizar and {Victoir}, Nicolas},
        title = "{Second order backward stochastic differential equations and fully non-linear parabolic PDEs}",
      journal = {arXiv Mathematics e-prints},
         year = 2005,
        month = sep,
          eid = {math/0509295},
        pages = {math/0509295},
          doi = {10.48550/arXiv.math/0509295},
archivePrefix = {arXiv},
       eprint = {math/0509295},
 primaryClass = {math.PR},
       adsurl = {https://ui.adsabs.harvard.edu/abs/2005math......9295C}
}

@article{Kingma2014AdamAM,
  title={Adam: A Method for Stochastic Optimization},
  author={Diederik P. Kingma and Jimmy Ba},
  journal={CoRR},
  year={2014},
  volume={abs/1412.6980},
  url={https://api.semanticscholar.org/CorpusID:6628106}
}

@inproceedings{he2016deep,
  author    = {He, Kaiming and Zhang, Xiangyu and Ren, Shaoqing and Sun, Jian},
  title     = {Deep residual learning for image recognition},
  booktitle = {IEEE Conference on Computer Vision and Pattern Recognition (CVPR)},
  year      = {2016}
}

@ARTICLE{2019arXiv191011623G,
       author = {{G{\"u}ler}, Batuhan and {Laignelet}, Alexis and {Parpas}, Panos},
        title = "{Towards Robust and Stable Deep Learning Algorithms for Forward Backward Stochastic Differential Equations}",
      journal = {arXiv e-prints},
         year = 2019,
        month = oct,
          eid = {arXiv:1910.11623},
        pages = {arXiv:1910.11623},
          doi = {10.48550/arXiv.1910.11623},
archivePrefix = {arXiv},
       eprint = {1910.11623},
 primaryClass = {stat.ML},
       adsurl = {https://ui.adsabs.harvard.edu/abs/2019arXiv191011623G}
}

@article{55960847-f1d7-372c-8587-c61145595486,
 ISSN = {00221082, 15406261},
 URL = {http://www.jstor.org/stable/2328113},
 author = {Hayne E. Leland},
 journal = {The Journal of Finance},
 number = {5},
 pages = {1283--1301},
 publisher = {[American Finance Association, Wiley]},
 title = {Option Pricing and Replication with Transactions Costs},
 urldate = {2026-07-31},
 volume = {40},
 year = {1985}
}

@ARTICLE{Barles1998-sz,
  title     = "Option pricing with transaction costs and a nonlinear
               {Black-Scholes} equation",
  author    = "Barles, Guy and Soner, Halil Mete",
  journal   = "Finance Stoch.",
  publisher = "Springer Science and Business Media LLC",
  volume    =  2,
  number    =  4,
  pages     = "369--397",
  month     =  aug,
  year      =  1998
}

@ARTICLE{Lepinette2012-ae,
  title     = "Modified Leland's strategy for a constant transaction costs rate",
  author    = "Lepinette, Emmanuel",
  journal   = "Math. Financ.",
  publisher = "Wiley",
  volume    =  22,
  number    =  4,
  pages     = "741--752",
  month     =  oct,
  year      =  2012,
  copyright = "http://onlinelibrary.wiley.com/termsAndConditions\#vor",
  language  = "en"
}

@article{92dffd60-7698-3838-9cfd-f3a5e322ac6b,
 ISSN = {00029947, 10886850},
 URL = {http://www.jstor.org/stable/1990512},
 author = {M. Kac},
 journal = {Transactions of the American Mathematical Society},
 number = {1},
 pages = {1--13},
 publisher = {American Mathematical Society},
 title = {On Distributions of Certain Wiener Functionals},
 urldate = {2026-08-13},
 volume = {65},
 year = {1949}
}

@inproceedings{Pardoux1992BackwardSD,
  title={Backward stochastic differential equations and quasilinear parabolic partial differential equations},
  author={{\'E}tienne Pardoux and Shige Peng},
  year={1992},
  url={https://api.semanticscholar.org/CorpusID:119715750}
}

@article{0b9b8115-a8b8-3422-8e1c-a62077de6621,
 ISSN = {00223808, 1537534X},
 URL = {http://www.jstor.org/stable/1831029},
 author = {Fischer Black and Myron Scholes},
 journal = {Journal of Political Economy},
 number = {3},
 pages = {637--654},
 publisher = {The University of Chicago Press},
 title = {The Pricing of Options and Corporate Liabilities},
 urldate = {2026-08-14},
 volume = {81},
 year = {1973}
}

@article{b5405564-0209-3283-bc5e-72bd72cfbbed,
 ISSN = {00058556, 23255323},
 URL = {http://www.jstor.org/stable/3003143},
 author = {Robert C. Merton},
 journal = {The Bell Journal of Economics and Management Science},
 number = {1},
 pages = {141--183},
 publisher = {[Wiley, RAND Corporation]},
 title = {Theory of Rational Option Pricing},
 urldate = {2026-08-14},
 volume = {4},
 year = {1973}
}

@ARTICLE{Saeed2026-cg,
  title     = "Option pricing under stochastic interest rates and transaction
               costs using physics-informed neural networks",
  author    = "Saeed, Tareq and Hussain, Javed and Rafique, Sumaiya",
  journal   = "Examples and Counterexamples",
  publisher = "Elsevier BV",
  volume    =  10,
  number    =  100227,
  pages     = "100227",
  month     =  dec,
  year      =  2026,
  copyright = "http://creativecommons.org/licenses/by-nc-nd/4.0/",
  language  = "en"
}

@ARTICLE{Herdegen2021-nf,
  title     = "Equilibrium asset pricing with transaction costs",
  author    = "Herdegen, Martin and Muhle-Karbe, Johannes and Possama{\"\i},
               Dylan",
  journal   = "Finance Stoch.",
  publisher = "Springer Science and Business Media LLC",
  volume    =  25,
  number    =  2,
  pages     = "231--275",
  month     =  apr,
  year      =  2021,
  copyright = "https://creativecommons.org/licenses/by/4.0",
  language  = "en"
}

@inproceedings{smith2017cyclical,
  title={Cyclical learning rates for training neural networks},
  author={Smith, Leslie N},
  booktitle={2017 IEEE winter conference on applications of computer vision (WACV)},
  pages={464--472},
  year={2017},
  organization={IEEE}
}

@Article{info11020108,
AUTHOR = {Howard, Jeremy and Gugger, Sylvain},
TITLE = {Fastai: A Layered API for Deep Learning},
JOURNAL = {Information},
VOLUME = {11},
YEAR = {2020},
NUMBER = {2},
ARTICLE-NUMBER = {108},
URL = {https://www.mdpi.com/2078-2489/11/2/108},
ISSN = {2078-2489},
DOI = {10.3390/info11020108}
}

@ARTICLE{2019JNS....29.1563B,
       author = {{Beck}, Christian and {E}, Weinan and {Jentzen}, Arnulf},
        title = "{Machine Learning Approximation Algorithms for High-Dimensional Fully Nonlinear Partial Differential Equations and Second-order Backward Stochastic Differential Equations}",
      journal = {Journal of NonLinear Science},
         year = 2019,
        month = aug,
       volume = {29},
       number = {4},
        pages = {1563-1619},
          doi = {10.1007/s00332-018-9525-3},
archivePrefix = {arXiv},
       eprint = {1709.05963},
 primaryClass = {cs.NA},
       adsurl = {https://ui.adsabs.harvard.edu/abs/2019JNS....29.1563B}
}

@ARTICLE{2019arXiv190800412P,
       author = {{Pham}, Huyen and {Warin}, Xavier and {Germain}, Maximilien},
        title = "{Neural networks-based backward scheme for fully nonlinear PDEs}",
      journal = {arXiv e-prints},
         year = 2019,
        month = jul,
          eid = {arXiv:1908.00412},
        pages = {arXiv:1908.00412},
          doi = {10.48550/arXiv.1908.00412},
archivePrefix = {arXiv},
       eprint = {1908.00412},
 primaryClass = {math.OC},
       adsurl = {https://ui.adsabs.harvard.edu/abs/2019arXiv190800412P}
}

@ARTICLE{2026arXiv260716730Z,
       author = {{Zhao}, Zhenhua and {Long}, Jihao},
        title = "{A Deep Second-Order Stochastic Residual Method for Fully Nonlinear Parabolic PDEs}",
      journal = {arXiv e-prints},
         year = 2026,
        month = jul,
          eid = {arXiv:2607.16730},
        pages = {arXiv:2607.16730},
          doi = {10.48550/arXiv.2607.16730},
archivePrefix = {arXiv},
       eprint = {2607.16730},
 primaryClass = {cs.NA},
       adsurl = {https://ui.adsabs.harvard.edu/abs/2026arXiv260716730Z}
}

@article{Marco01121994,
author = {Avellaneda, Marco and Par{á}s, Antonio},
title = {Dynamic hedging portfolios for derivative securities in the presence of large transaction costs},
journal = {Applied Mathematical Finance},
volume = {1},
number = {2},
pages = {165--194},
year = {1994},
publisher = {Routledge},
doi = {10.1080/13504869400000010},


URL = { 
    
        https://doi.org/10.1080/13504869400000010
    
    

},
eprint = { 
    
        https://doi.org/10.1080/13504869400000010
    
    

}

}

@misc{sitzmann2020implicitneuralrepresentationsperiodic,
      title={Implicit Neural Representations with Periodic Activation Functions}, 
      author={Vincent Sitzmann and Julien N. P. Martel and Alexander W. Bergman and David B. Lindell and Gordon Wetzstein},
      year={2020},
      eprint={2006.09661},
      archivePrefix={arXiv},
      primaryClass={cs.CV},
      url={https://arxiv.org/abs/2006.09661}, 
}
\bibliographystyle{icml2025}

\newpage
\appendix
\onecolumn
\section{Network and training details}
\label{app:hyper}
In this Appendix we provide further details of the network and its training.

The input to the network is normalised underlying prices and time $(t,\vec{X_t}/\vec{X}_0)\in\mathbb R^{1+D}$. This is fed into a linear layer of width 256 with sine activation, which is then followed by three residual blocks used in the construction of a ResNet each with width 256. These map the inputs $h\mapsto h+\sin(Wh+b)$, where $W$ and $b$ are the weights and biases. The weights were initialised by the Xavier uniform scheme. The output layer is a linear layer with $1+D(D+1)/2$ outputs corresponding to the solution $u^\theta_t$ and the degrees of freedom in $\mat\Gamma^\theta_n$. There is an option to scale down the outputs of $\mat\Gamma^\theta_n$ by $\xi_i\xi_j$ (the initial values of $X_i$) so that the output learns the cash gamma at the initial prices. This is used only for the experiments at $D=100$.

The training proceeds by the simulation of $M=100$ fresh underlying paths with $N=50$ time steps (apart from $t=0$) using the Euler--Maruyama scheme in \cref{eq:LelandFwd}. The loss function in \cref{eq:loss} is evaluated and recorded. The training is designed such that any gradients associated with $\mat\Gamma^\theta$ are not passed to $\varphi_n$, the drift of the backward equation, since this tends to bias the learnt solution. The separation of the training is similar in spirit to the two separate neural network approaches traditionally employed in the 2BSDE literature, but our approach uses a single neural network instead. The Adam optimiser is chosen with default parameters and the learning rate $\eta^\ast$ is the ``valley'' output of the \texttt{fastai} learning rate range test, which tests learning rates from $10^{-6}$ up to $10^{-1}$ in 1000 iterations.

The evaluation of the network uses 100 newly generated test paths with $N=50$. Any quantity labelled by ``final'' refers to the network evaluating it after its final iteration of learning, whereas ``mean'' refers to the mean value ($\pm$ standard deviation) of values logged every 100 iterations over the last 2000 iterations.

A summary of the hyperparameters used in this work are presented in \cref{tab:hyper}.
\begin{table}[h]
\caption{Summary of hyperparameters for each experiment in this work. $^\dagger$Trained without stopping the gradients of $\mat\Gamma^\theta_n$ through $\varphi_n$.}
\label{tab:hyper}
\vskip 0.1in
\begin{center}
\begin{small}
\begin{tabular*}{\textwidth}{@{\extracolsep{\fill}}lrccrcrrr@{}}
\toprule
Product & $D$ & $\Gamma$-head & $\Gamma$ rescaled & Parameters & $\eta^*$ & Iterations (stage 1 + 2) & Time & s/100 it. \\
\midrule
Call                         & 10  & yes & no  & 214\,840   & $1.51\times10^{-4}$ & 21\,600 + 4\,300 & 2.5 h  & 35 \\
Butterfly                    & 10  & yes & no  & 214\,840   & $1.43\times10^{-4}$ & 16\,800 + 2\,900 & 2.0 h  & 36 \\
Butterfly (no $\Gamma$-head) & 10  & no  & --  & 214\,840   & $2.16\times10^{-4}$ & 24\,400 + 3\,200 & 15.0 h & 195 \\ 
Call                         & 100 & yes & yes & 1\,521\,595 & $2.95\times10^{-4}$ & 21\,500 + 6\,100 & 11.4 h & 149 \\
Butterfly                    & 100 & yes & yes & 1\,521\,595 & $3.72\times10^{-4}$ & 24\,900 + 7\,300 & 12.9 h & 144 \\
\midrule
Butterfly$^\dagger$          & 10  & yes & no  & 214\,840   & $1.19\times10^{-4}$ & 22\,900 + 3\,000 & 2.9 h  & 40 \\
\bottomrule
\end{tabular*}
\end{small}
\end{center}
\vskip -0.1in
\end{table}
\section{Geometric baskets}
\label{sec:GeomBasket}
The geometric basket is a heavily utilised concept in this work and the aim of this Appendix is to provide some additional context about it.

In this work we consider $D$ lognormally distributed assets in the risk-neutral framework without correlation between Brownian motions. They obey the following SDE \begin{equation}
    d\vec X_t=r\vec X_tdt +\mat\sigma(\vec X_t)\cdot d\vec W_t
\end{equation}where $\mat\sigma(\vec X_t)=\operatorname{diag}(\sigma_iX_i)$. Solving this SDE gives the following, 
\begin{equation}
    \ln X_i(t)=\ln X_i(0)+(r-\tfrac12\sigma_i^2)t+\sigma_iW_i(t)
\end{equation}
The geometric basket is defined as 
\begin{equation}
    G_t=\prod_{i=1}^DX_i(t)^{w_i}
    \end{equation}
which can also be expressed as \begin{equation}
    \ln G_t = \ln G_0 +\left(r-\frac{1}{2}\sum_iw_i\sigma_i^2\right)t + \sum_iw_i\sigma_iW_i(t)\label{eq:lnGt}
\end{equation}
This reduces the $D$-dimensional problem into a one-dimensional problem since the basket is a scalar quantity which obeys its own geometric Brownian motion given by the SDE (with a new Brownian motion) \begin{equation}
    dG_t = (r-q_G)G_tdt+\sigma_G G_tdW_t,\quad q_G= \frac12\left(\sum_iw_i\sigma_i^2-\sigma_G^2\right),\;\sigma_G^2=\sum_iw_i^2\sigma_i^2\label{eq:GSDE}
\end{equation}
Since each asset has drift $r$, arithmetic basket also has drift $r$. However, the drift of a geometric basket should be less than $r$ by the AM--GM inequality and so the difference $q_G\geq0$ appears like a dividend in the SDE.

Therefore the price of a product on a geometric basket of $D$ underlyings are equivalently solved by the one-dimensional version using the basket SDE above. Or through the following Black--Scholes equation for the derivative price $u(t,\vec X_t)=V(t,G_t)$ \begin{equation}
    \frac{\partial V}{\partial t} + \frac{1}2\sigma_G^2G^2\frac{\partial^2V}{\partial G^2}+(r-q_G)G\frac{\partial V}{\partial G}-rV=0
\end{equation}

\subsection{The Leland term for geometric baskets}
\label{sec:LelandGeom}
Vital to the Leland transaction model is the sign of the diagonal components of the products curvature $\Gamma_{ii}$, which can be computed straight forwardly using the basket description as \begin{equation}
    \Gamma_{ii}=\frac{\partial^2}{\partial X_i^2}V(t,G_t)=\frac{w_i^2}{X_i^2}\left[G^2V_{GG}+\left(1-\frac{1}{w_i}\right)GV_G\right]=\frac{w_i^2}{X_i^2}\mathcal Q_i
\end{equation}
One then notices that \begin{equation}
    \frac{1}2\sigma_G^2G^2\frac{\partial^2V}{\partial G^2}-q_GG\frac{\partial V}{\partial G}=\frac{1}2\sum_i\sigma_i^2w_i^2\mathcal Q_i
\end{equation}
and so incorporating the Leland term gives 
\begin{equation}
    \frac{\partial V}{\partial t} +rG\frac{\partial V}{\partial G}+\frac{1}2\sum_i\sigma_i^2w_i^2\mathcal Q_i(1+\operatorname{Le}_i\operatorname{sign}\mathcal Q_i)-rV=0
\end{equation}

When setting equal $\sigma_i$ and $w_i=1/D$, above formulae simplify. In particular $\mathcal Q_i=\mathcal Q$ where \begin{equation}
    \mathcal Q=G^2\frac{\partial^2 V}{\partial G^2}-(D-1)G\frac{\partial V}{\partial G}\label{eq:Q}
\end{equation}and the one-dimensional PDE becomes 
\begin{equation}
    \frac{\partial V}{\partial t} +rG\frac{\partial V}{\partial G}+\frac{1}2\sigma_G^2\mathcal Q(1+\operatorname{Le}\operatorname{sign}\mathcal Q)-rV=0 \label{eq:BS1D}
\end{equation}

\section{Finite-difference reference solver}
\label{app:fd}
The main reference for the PDE solutions comes from a finite difference solver. In this Appendix we specify some details about the scheme and grid size for this work.

We assume that $w_i=1/D$, $\sigma_i=\sigma$ and $\mathrm{Le}_i=\mathrm{Le}$, and $\mathcal{Q}_i=\mathcal Q$ as in \cref{eq:Q}. The one-dimensional Black--Scholes equation takes the form in \cref{eq:BS1D} with a solution $V(t,G)$ with $G$ as the geometric basket level. It is convenient to define the variables $y=\ln G$ and $\tau=T-t$ to re-express the PDE as

\begin{equation}
\frac{\partial V}{\partial \tau} = r\frac{\partial V}{\partial y}+\frac{\sigma^2}{2D}\bigl(\mathcal{Q}+\mathrm{Le}\,|\mathcal{Q}|\bigr)-rV,
\qquad \mathcal{Q}=\frac{\partial^2V}{\partial y^2}-D\frac{\partial V}{\partial y}
\label{eq:fd_pde}
\end{equation}

\cref{eq:fd_pde} is integrated with the explicit Euler method in $\tau$ and second-order central differences for the derivatives. The nonlinearity in $|\mathcal{Q}|$ is handled by evaluating it pointwise at the previous time step.

Choosing appropriate step sizes $\Delta y$ and $\Delta\tau$ requires care and will be fixed by a von Neumann analysis to ensure stability of the solution. For convenience we define $\kappa_\pm=\sigma^2(1\pm\mathrm{Le})/(2D)$ as the coefficient of the $V_{yy}$ term in the PDE and $a_\pm=r-\tfrac12\sigma^2(1\pm\mathrm{Le})$ for the coefficient of the $V_y$. At any given time step, say node $j$, these coefficients are fixed which we refer to without the subscript. The explicit finite difference equation is written as
\begin{equation}
V_j^{n+1}=V_j^{n}+\Delta\tau\left[a\,\frac{V_{j+1}^n-V_{j-1}^n}{2\Delta y}+\kappa\,\frac{V_{j+1}^n-2V_j^n+V_{j-1}^n}{\Delta y^2}-rV_j^n\right].
\end{equation}
With the definition of the mesh ratio $\lambda=\kappa\Delta\tau/\Delta y^2$, the Courant number $\mu=a\Delta\tau/\Delta y$ and collecting terms gives
\begin{equation}
V_j^{n+1}=\bigl(\lambda+\tfrac{\mu}{2}\bigr)V_{j+1}^n+\bigl(1-2\lambda-r\Delta\tau\bigr)V_j^n+\bigl(\lambda-\tfrac{\mu}{2}\bigr)V_{j-1}^n .
\label{eq:fd_update}
\end{equation}
The von Neumann analysis proceeds with the substitution of a Fourier mode $V_j^n=g^ne^{ij\theta}$ restricting $\theta\in[-\pi,\pi]$ into the \cref{eq:fd_update} to extract the amplification factor\begin{equation}
g(\theta)=1-r\Delta\tau-4\lambda\sin^2\tfrac{\theta}{2}+i\,\mu\sin\theta
\end{equation}
For the stability analysis the $-r\Delta\tau$ term can be dropped since it damps the solution and is of higher order so it does not contribute to the stability of the solution. Writing $s=\sin^2(\theta/2)\in[0,1]$ and using $\sin^2\theta=4s(1-s)$,
\begin{equation}
|g|^2-1=(1-4\lambda s)^2+4\mu^2 s(1-s)-1=4s\Bigl[(\mu^2-2\lambda)+(4\lambda^2-\mu^2)\,s\Bigr]
\end{equation}
The von Neumann stability condition for stability requires $|g|\leq1\,\forall\theta$. This essentially requires the term in brackets to be non-positive giving two conditions
\begin{equation}
\mu^2\le2\lambda\ \ \Leftrightarrow\ \ \Delta\tau\le\frac{2\kappa}{a^2},
\qquad\qquad
\lambda\le\tfrac12\ \ \Leftrightarrow\ \ \Delta\tau\le\frac{\Delta y^2}{2\kappa}.
\label{eq:vn_conditions}
\end{equation}
and puts the following condition on $\Delta \tau$ and $\Delta y$\begin{equation}
    \Delta \tau \leq \operatorname{min}\left(\frac{2\kappa}{a^2},\frac{\Delta y^2}{2\kappa}\right),\qquad \Delta y\leq\frac{2\kappa}{|a|}
\end{equation}which is saturated at $\Delta y=2\kappa/|a|$. The bound on $\Delta y$ follows from monotonicity.

These are bounds for a given time step where $\kappa$ and $a$ have been fixed, however the lowest possible upper bound on $\Delta y$ is\begin{equation}
\Delta y\le\min_\pm\frac{2\kappa_\pm}{|a_\pm|}=4.132\times10^{-3}
\end{equation}
when choosing the ``$-$ '' values and $D=100$. Once $\Delta y$ is fixed one can fix $\Delta t$ as
\begin{equation}
\Delta\tau\le\min\!\left(\frac{2\kappa_+}{a_+^2},\frac{2\kappa_-}{a_-^2},\ \frac{\Delta y^2}{2\kappa_+}\right)
\end{equation}In the examples studied here one has $2\kappa_+/a_+^2 >2\kappa_-/a_-^2$.

The grid used here has $N_y=4001$ equally spaced points on $y\in[\ln G_0-1.5,\ \ln G_0+1.5]$ or equivalently $G\in[22.3,\,448.2]$. The separation in $y$ is $\Delta y=7.5\times10^{-4}$ which is below the largest possible value for $D=100$. The boundaries lie more than 20 standard deviations of $\ln G_T$ from $\ln G_0$ for every case considered. This choice also ensures that $\Delta y^2/2\kappa_+<2\kappa_-/a_-^2$ even for $D=100$.

Replacing $\tfrac12$ by a safety factor of $0.4$ gives the stability limit $\Delta\tau_{\mathrm{stab}}=0.4\,\Delta y^2/\kappa_+$. Having dealt with the stability of the solution it is time to address the accuracy of the solution. The leading truncation error of the scheme removes a fraction $a_\pm^2\Delta\tau/(2\kappa_\pm)$ of the diffusion. Based on some scaling arguments, since $\kappa_\pm\propto1/D$ while  $\Delta\tau_{\mathrm{stab}}\propto D$, this fraction grows like $D^2$. Taking just the stability on its own gives  a proportion of about $2\times10^{-4}$ at $D=10$ but up to $0.02$ at $D=100$. We therefore also impose the accuracy limit
\begin{equation}
\Delta\tau_{\mathrm{stab}}=\frac{0.4\,\Delta y^2}{\kappa_+},\qquad
\Delta\tau_{\mathrm{acc}}=\frac{2\varepsilon\,\kappa_-}{a_-^2},\qquad \varepsilon=10^{-3},
\end{equation}
and take
\begin{equation}
N_\tau=\max\!\left(\left\lceil\frac{T}{\min(\Delta\tau_{\mathrm{stab}},\,\Delta\tau_{\mathrm{acc}})}\right\rceil+1,\ 400\right),\qquad
\Delta\tau=\frac{T}{N_\tau}.
\end{equation}
At $D=10$ the stability limit is smaller and so $N_\tau=10\,169$, whereas at $D=100$ the accuracy limit is smaller and $N_\tau=20\,058$.

\end{document}